\documentclass[
submission
]{dmtcs-episciences}

\usepackage[utf8]{inputenc}
\usepackage{subfigure}

\usepackage[round]{natbib}

\author{Samvel Kh. Darbinyan}
\title[A new sufficient condition 
 for a digraph to be   Hamiltonian$-$A proof of Manoussakis conjecture ]{A new sufficient condition 
 for a digraph to be   Hamiltonian$-$A proof of Manoussakis conjecture }
\affiliation{  Yerevan, Armenia\\Institute for Informatics and Automation Problems of NAS RA}
\received{2020-02-11}
\revised{2020-07-21,2020-11-05}
\accepted{2020-11-05}

\begin{document}
\publicationdetails{22}{2020}{4}{12}{6086}
\maketitle
\begin{abstract}
   Y. Manoussakis (J. Graph Theory 16, 1992, 51-59) proposed the following conjecture.

\noindent\textbf{Conjecture}. {\it Let $D$ be a 2-strongly connected  digraph of order $n$ such that
for all distinct pairs of non-adjacent vertices $x$, $y$ and $w$, $z$, we have $d(x)+d(y)+d(w)+d(z)\geq 4n-3$. Then $D$ is Hamiltonian.}

In this paper, we confirm this  conjecture. Moreover, we prove  that if a  digraph $D$ satisfies the conditions of this conjecture and   has a pair of non-adjacent vertices  $\{x,y\}$ such that  $d(x)+d(y)\leq 2n-4$,  then $D$ contains cycles of all lengths $3, 4, \ldots , n$.
\end{abstract}

\section{Introduction}

 In this paper, we consider finite digraphs (directed graphs)  without loops and multiple arcs. Every cycle and path are assumed simple and directed; its {\it length} is the number of its arcs.                                                                                                
A digraph $D$ is {\it Hamiltonian} 
 if it contains a cycle passing through all the vertices of $D$. There are many conditions that guarantee that  a digraph is Hamiltonian (see, e.g., \cite{[1]}, \cite{[3]}, \cite{[17]}, \cite{[19]}, \cite{[20]}).  \cite{[19]}  the following theorem was proved.\\

  \noindent\textbf{Theorem 1.1} (\cite{[19]}). {\it Let $D$ be a strong  digraph of order $n\geq 4$. Suppose that $D$ satisfies the following condition for every triple $x,y,z\in V(D)$ such that $x$ and $y$ are non-adjacent: If there is no arc from $x$ to $z$, then
  $d(x)+d(y)+d^+(x)+d^-(z)\geq 3n-2$.  If there is no arc from $z$ to $x$, then
  $d(x)+d(y)+d^-(x)+d^+(z)\geq 3n-2$. Then $D$ is Hamiltonian.}\\

  \noindent\textbf{Definition 1.2.}  {\it Let $D$ be a digraph of order $n$. We say that $D$ satisfies condition $(M)$ when $d(x)+d(y)+d(w)+d(z)\geq 4n-3$ for all distinct pairs of non-adjacent vertices $x,y$ and $w,z$.}\\

  \cite{[19]} proposed the following conjecture. This conjecture is an extension of Theorem 1.1. 

  \noindent\textbf{Conjecture 1.3} (\cite{[19]}). {\it Let $D$ be a 2-strong  digraph of order $n$ such that
for all distinct pairs of non-adjacent vertices $x$, $y$ and $w$, $z$ we have $d(x)+d(y)+d(w)+d(z)\geq 4n-3$. Then $D$ is Hamiltonian.}\\

 \cite{[19]}
gave an example, which  showed that if this  conjecture is true, then the minimum degree condition is sharp. Notice that another examples can be found  in a paper by \cite{[7]}, where  
 for any two integers $k\geq 2$ and $m\geq 1$,  the author  constructed a family of  $k$-strong  digraphs of order $4k+m$ with  minimum degree  $4k+m-1$, which are not Hamiltonian. This result improves a conjecture of Thomassen (see \cite{[3]} Conjecture 1.4.1:  Every 2-strong $(n-1)$-regular digraph of order $n$, except $D_5$ and $D_7$, is Hamiltonian). Moreover, when $m=1$,  then from these digraphs we can obtain $k$-strong non-Hamiltonian digraphs of order $n=4k+1$ with  minimum degree equal to $n-1$ and the minimal semi-degrees equal to $(n-3)/2$. Thus, if in Conjecture 1.3 we replace $4n-3$ with $4n-4$, then for every $n$ there are many digraphs of order $n$ with high connectivity and  high semi-degrees, for which Conjecture 1.3 is not true. 
 
 The {\it cycle factor} in a digraph $D$ is a collection of pairwise vertex disjoint cycles $C_1, C_2, \dots , C_l$ such that $\bigcup ^l_{i=1}V(C_i)=V(D)$. It is clear that the existence of a cycle factor in a digraph $D$ is a necessary condition for a digraph to be Hamiltonian. The following theorem   gives a necessary and sufficient condition for the existence of a cycle factor in a digraph.\\
 
 \noindent\textbf{Theorem 1.4} (\cite{[25]}). {\it Let $D$ be a digraph. Then $D$ has a cycle factor if and only if
$V(D)$ cannot be partitioned into subsets $Y$, $Z$, $R_1$, $R_2$ such that $A(Y\rightarrow R_1)=A(R_2\rightarrow R_1\cup Y)=\emptyset$, $|Y|>|Z|$ and $Y$ is an independent set.}\\

Using  theorem Theorem 1.4,  it is not difficult to construct  2-strong digraphs satisfying the condition that 
 $d(x)+d(y)+d(w)+d(z)\geq 4n-4$ for every distinct pairs $\{x,y\}$, $\{w,z\}$ of non-adjacent vertices, but these digraphs do not even contain a cycle factor. 

Thomassen  suggested (see \cite{[3]}) the following two conjectures:\\

\noindent\textbf{1}. Conjecture 1.6.7.  {\it Every 3-strong digraph of order $n$ and with minimum degree at least $n+1$
 is strongly Hamiltonian-connected}.\\

\noindent\textbf{2}. Conjecture 1.6.8.  {\it Let $D$ be a 4-strong digraph of order $n$ such that the sum of the degrees of any pair of non-adjacent vertices  is at least $2n+1$. Then $D$ is strongly Hamiltonian-connected}.\\

Investigating these conjectures,  \cite{[8]} disproved the first conjecture (proving that for every integer $n\geq 9$ there exists a 3-strong non-strongly Hamiltonian-connected digraph of order $n$ with the minimum degree at least $n+1$)  and for the second proved 
 the following two theorems.\\
 
\noindent\textbf{Theorem 1.5.}  {\it Any  $k$-strong ($k\geq 1$)  digraph $D$ of order $n\geq 8$ satisfying the condition that the sum of degrees of any pair of non-adjacent vertices $x, y\in V(D)\setminus \{z\}$  at least $2n-1$, where $z$ is some vertex in $V(D)$, is Hamiltonian if and only if any $(k+1)$-strong digraph of order $n+1$ satisfying the condition that the sum of degrees of any pair of non-adjacent vertices at least $2n+3$ is strongly Hamiltonian-connected}.\\

\noindent\textbf{Theorem 1.6.}  {\it Let $D$ be a strong  digraph of order $n\geq 3$. Suppose that $d(x)+d(y)\geq 2n-1$ for every pair of non-adjacent vertices $x, y\in V(D)\setminus \{z\}$, where $z$ is some vertex of $V(D)$. Then  $D$ contains a cycle of length  at least  $n-1$.}\\

It is easy to see that if a digraph $D$ satisfies the condition $(M)$, then   it contains at most one pair of non-adjacent vertices $x,y$ such that $d(x)+d(y)\leq 2n-2$. From this and Theorem 1.6,   the following corollary immediately  follows.\\

\noindent\textbf{Corollary 1.7.}  {\it Let $D$ be a strong digraph of order $n$ satisfying condition $(M)$.  Then $D$ contains a cycle of length at least $n-1$ (in particular, $D$ contains a Hamiltonian path).}\\

Corollary 1.7 was also  later  proved by  \cite{[22]}.

It is worth to note that   \cite{[9]}, \cite{[10]} (2015) and \cite{[11]}  studied  some properties in digraphs with the conditions of Theorem 1.1. They  obtained the following results (in first two results $D$ is a digraph of order $n$ satisfying the degree condition of Theorem 1.1). \\

(i) (\cite{[11]}). {\it If $D$ is strong, then it contains a cycle of length $n-1$ or   $D$ is isomorphic to the complete bipartite digraph $K^*_{n/2,n/2}$}.\\
 
 (ii) (\cite{[10]} (2015)). {\it If  $D$ is strong, then it contains a Hamiltonian path in which the initial vertex dominates the terminal vertex  or $D$ is isomorphic to one tournament of order  5}.\\

(iii) (\cite{[9]}). {\it Let $D$ be  a digraph of order $n$ and let $Y$ be a non-empty subset of $V(D)$. Suppose that for every triple  of the vertices $x,y,z\in Y$ such that $x$ and $y$ are non-adjacent:  If there is no arc from $x$ to $z$, then
  $d(x)+d(y)+d^+(x)+d^-(z)\geq 3n-2$.  If there is no arc from $z$ to $x$, then
  $d(x)+d(y)+d^-(x)+d^+(z)\geq 3n-2$. If  there is a path from $u$ to $v$ and a path from $v$ to $u$ in $D$ for every pair of distinct vertices  $u,v\in Y$, then $D$ has a cycle which contains at least $|Y|-1$  vertices of $Y$}.\\

The last result is best possible in some situations and gives an answer to a question posted by   \cite{[18]}.

\noindent\textbf{Theorem 1.8} (\cite{[20]}). {\it Let $D$ be a strong  digraph of order $n\geq 2$. If $d(x)+d(y)\geq 2n-1$ for all pairs of non-adjacent vertices $x, y$ in $D$, then $D$ is Hamiltonian.}\\

For a short proof of Theorem 1.8, see \cite{[4]}.
 \cite{[6]} characterized those digraphs which satisfy Meyniel's condition, but are not pancyclic.
Before stating the main result obtained  by  \cite{[6]}, we need to define a family of digraphs.\\

\noindent\textbf{Definition 1.9.} {\it For  integers $n$ and $m$, $(n+1)/2 < m\leq n-1$, let $\Phi_n^m$ denote the set of
 digraphs $D$, which satisfy the following conditions: (i) $V(D)= \{x_1,x_2,\ldots , x_n\}$; (ii) $x_nx_{n-1}\ldots x_2x_1x_n$ is a Hamiltonian cycle in $D$; (iii) for each $k$, $1\leq k\leq n-m+1$, the vertices $x_k$ and $x_{k+m-1}$ are not adjacent;
(iv) $x_jx_i\notin A(D)$ whenever $2\leq i+1< j\leq n$ and (v) the sum of degrees for any two distinct non-adjacent vertices is at least $2n-1$.}\\

\noindent\textbf{Theorem 1.10} (\cite{[5]}, \cite{[6]}). {\it Let $D$ be a strong digraph of order $n\geq 3$. Suppose that $d(x)+d(y)\geq 2n-1$ for all pairs of distinct non-adjacent vertices $x$, $y$ in $D$. Then either (a) $D$ is pancyclic or (b) $n$ is even and $D$ is isomorphic to one of $K^*_{n/2,n/2}$, $K^*_{n/2,n/2}\setminus \{e\}$, where $e$ is an arbitrary arc of $K^*_{n/2,n/2}$, or (c) $D\in \Phi^m_n $ $($in this case $D$ does not contain a cycle of length $m$}).\\

 Later on, Theorem 1.10  was also proved by  \cite{[2]}.  \cite{[13]} investigated the pancyclicity of digraphs with the condition $(M)$. Using Theorem 1.10 and  the Moser theorem for a strong tournament to be  pancyclic (see \cite{[16]}), we proved the following theorem.\\
 
 \noindent\textbf{Theorem 1.11} (\cite{[13]}). {\it Let $D$ be a 2-strong digraph of order $n\geq 6$ satisfying condition $(M)$. Suppose that there exists 
  a pair   of non-adjacent vertices $x$, $y$ in $D$ such that 
  $d(x)+d(y)\leq 2n-4$. Then $D$ contains cycles of all lengths $3, 4, \ldots , n-1$}.\\

  In this paper we confirm  Conjecture 1.3.\\
  
  \noindent\textbf{Theorem 1.12.} {\it Let $D$ be a 2-strong digraph of order $n\geq 3$ satisfying condition $(M)$. Then $D$ is Hamiltonian}.

  Theorem 1.12 also has the following immediate corollaries.\\
  
    \noindent\textbf{Corollary 1.13} (\cite{[24]}).
     {\it A digraph of order $n$ is Hamiltonian if, for any two vertices $x$ and $y$, either $x \rightarrow y$ or $d^+(x)+d^-(y)\geq n$}.\\
    
    \noindent\textbf{Corollary 1.14} (\cite{[21]}). {\it Let $D$ be a digraph of order $n\geq 2$. If for every vertex $x$,  $d^+(x)\geq n/2$ and $d^-(x)\geq n/2$, then $D$ is Hamiltonian}.\\
    
 Note that Corollary 1.14 immediately follows from well-known
 theorem of Ghouila-Houri \cite{[14]}.\\
   
    \noindent\textbf{Corollary 1.15} ( \cite{[23]}). {\it Let $G$ be a simple graph of order $n\geq 3$, in which the degree sum of any two non-adjacent vertices is at least $n$. Then $G$ is Hamiltonian}. \\

 As an immediate corollary of Theorems 1.12 and 1.11, we obtain the following theorem.\\
 
 \noindent\textbf{Theorem 1.16.} {\it Let $D$ be a 2-strong digraph of order $n\geq 6$ satisfying condition $(M)$. Suppose that $D$ contains 
  a pair   of non-adjacent vertices $x$, $y$ such that 
  $d(x)+d(y)\leq 2n-4$. Then $D$ is pancyclic}.\\
   
In view of Theorem 1.16, it is natural to set the following problem.\\

\noindent\textbf{Problem 1.17.} {\it Let $D$ be a 2-strong connected digraph of order $n$ satisfying condition $(M)$. Suppose that $\{x,y\}$ is a pair of non-adjacent vertices in $D$ such that $2n-3\leq d(x)+d(y)\leq 2n-2$. Whether $D$ is pancyclic?}

\section{Terminology and notation}

In this paper we consider finite digraphs without loops and multiple arcs. We shall assume that the reader is familiar with the standard  terminology on digraphs and refer to  the book \cite{[1]} for terminology and notations not defined here.
The vertex set and the arc set of a digraph $D$ are    denoted
  by $V(D)$  and   $A(D)$, respectively.  The {\it order} of $D$ is the number
  of its vertices. For any $x,y\in V(D)$, we also write $x\rightarrow y$ if $xy\in A(D)$. We use the notations $\overrightarrow{a} [x,y]=1$ if $xy\in A(D)$ and $\overrightarrow{a} [x,y]=0$ if $xy\notin A(D)$.
 If $xy\in A(D)$,  $y$ is an {\it out-neighbour} of $x$ and $x$ is an {\it in-neighbour} of $y$.
 If $x\rightarrow y$ and $y\rightarrow z$, we write $x\rightarrow y\rightarrow z$. Two distinct vertices $x$ and $y$ are {\it adjacent} if $xy\in A(D)$ or $yx\in A(D)$ (or both).
If there is no arc from $x$ to $y$, we shall use the notation $xy\notin A(D)$.

 We let $N^+(x)$, $N^-(x)$ denote the set of  {\it out-neighbours}, respectively the set  of {\it in-neighbours} of a vertex $x$ in a digraph $D$.  If $A\subseteq V(D)$, then $N^+(x,A)= A \cap N^+(x)$ and $N^-(x,A)=A\cap N^-(x)$.
The {\it out-degree} of $x$ is $d^+(x)=|N^+(x)|$ and $d^-(x)=|N^-(x)|$ is the {\it in-degree} of $x$. Similarly, $d^+(x,A)=|N^+(x,A)|$ and $d^-(x,A)=|N^-(x,A)|$.
The {\it degree} of the vertex $x$ in $D$ is defined as $d(x)=d^+(x)+d^-(x)$ (similarly, $d(x,A)=d^+(x,A)+d^-(x,A)$). The subdigraph of $D$ induced by a subset $A$ of $V(D)$ is denoted by $D\langle A\rangle$.
 If $z$ is a vertex of a digraph $D$, then the subdigraph $D\langle V(D)\setminus \{z\}\rangle$   is denoted by $D-z$.  

For integers $a$ and $b$, $a\leq b$, let $[a,b]$  denote the set of
all integers, which are not less than $a$ and are not greater than
$b$.

The path (respectively, the cycle) consisting of the distinct vertices $x_1,x_2,\ldots ,x_m$ ($m\geq 2 $) and the arcs $x_ix_{i+1}$, $i\in [1,m-1]$  (respectively, $x_ix_{i+1}$, $i\in [1,m-1]$, and $x_mx_1$), is denoted  by $x_1x_2\cdots x_m$ (respectively, $x_1x_2\cdots x_mx_1$).
We say that $x_1x_2\cdots x_m$ is a path from $x_1$ to $x_m$ or is an $(x_1,x_m)$-{\it path}. Let $x$ and $y$ be two distinct vertices of a digraph $D$. Cycle that passing through $x$ and $y$ in $D$, we denote by $C(x,y)$.

A cycle (respectively, a path) that contains  all the vertices of $D$, is a  {\it Hamiltonian cycle} (respectively, is a {\it Hamiltonian path}).
A digraph is {\it Hamiltonian} if it contains a Hamiltonian cycle. A digraph $D$ is {\it strongly Hamiltonian-connected} if, for every ordered pair $\{x,y\} $ of distinct vertices of $D$ there is a Hamiltonian path from $x$ to $y$. A digraph $D$ of order $n\geq 3$ is {\it pancyclic} if it contains cycles of all lengths $m$, $3\leq m\leq n$.
  For a cycle  $C=x_1x_2\cdots x_kx_1$ of length $k$, the subscripts considered modulo $k$, i.e., $x_i=x_s$ for every $s$ and $i$ such that  $i\equiv s\, (\hbox {mod} \,k)$.
If $P$ is a path containing a subpath from $x$ to $y$,  we let $P[x,y]$ denote that subpath. Similarly, if $C$ is a cycle containing vertices $x$ and $y$, $C[x,y]$ denotes the subpath of $C$ from $x$ to $y$. 
If $j<i$, then $\{x_i, \ldots , x_j\}=\emptyset$.

A digraph $D$ is {\it strongly connected} (or just {\it strong}), if there exists a path from $x$ to $y$ and a path from $y$ to $x$ for every pair of distinct vertices $x,y$. A digraph $D$ is $k$-{\it strongly  connected} (or  $k$-{\it strong}), where $k\geq 1$, if $|V(D)|\geq k+1$ and $D\langle V(D)\setminus A\rangle$ is strongly connected for any subset $A\subset V(D)$ of at most $k-1$ vertices.

For a pair of disjoint subsets $A$ and $B$ of $V(D)$, we define
$A(A\rightarrow B)=\{xy\in A(D)\,|\, x\in A, y\in B\}$ and
$A(A,B)=A(A\rightarrow B)\cup A(B\rightarrow A)$.

\section{Auxiliary known results}

\noindent\textbf{Lemma 3.1}  
 (\cite{[15]}).  {\it Let $D$ be a digraph of order $n\geq 3$ containing a cycle $C$ of length $m$, $m\in [2,n-1]$. Let $x$ be a vertex not contained in this cycle. If $d(x,V(C))\geq m+1$, then $D$ contains a cycle of length $k$ for all $k\in [2,m+1]$ }.\\

It is not difficult to prove the following lemma.\\

\noindent\textbf{Lemma 3.2.} {\it Let  $D$ be a  digraph of order $n$. Assume that $xy\notin A(D)$ and the vertices $x$, $y$ in $D$  satisfy the degree condition  $d^+(x)+d^-(y)\geq n-2+k$, where $k\geq 1$. Then $D$ contains at least $k$ internally disjoint $(x,y)$-paths of length two.}\\
 
The following results were proved by \cite{[13]} and its preliminary version  presented at Emil Artin International Conference (\cite{[12]}). \\

\noindent\textbf{Theorem 3.3.}  {\it  Let $D$ be a 2-strong digraph of order $n\geq 3$ satisfying condition $(M)$. Suppose that $\{x,y\}$ is a pair of non-adjacent vertices in $V(D)$ such that
$d(x)+d(y)\leq 2n-2$. Then $D$ is Hamiltonian if and only if
 $D$ contains a cycle  through the vertices $x$ and $y$}.\\

\noindent\textbf{Theorem 3.4.}  {\it Let $D$ be a 2-strong digraph of order $n\geq 3$. Suppose that $D$ contains at most one   pair of non-adjacent vertices. Then $D$ is Hamiltonian.}\\

\noindent\textbf{Remark.}  {\it There is a strong non-Hamiltonian digraph of order $n\geq 5$, which is not 2-strong and has exactly one pair of non-adjacent vertices.}\\

Using Lemma 3.2, it is not difficult to prove the following lemma.\\

\noindent\textbf{Lemma 3.5.}   {\it Let $D$ be a 2-strong digraph of order $n\geq 3$ and  let $u$, $v$ be two distinct vertices in $V(D)$. If $D$ contains no cycle  through $u$ and $v$, then $u$, $v$ are not adjacent and  there is no path of length two between them. In particular, $d(u)+d(v)\leq 2n-4$}. \\

\noindent\textbf{Theorem 3.6.}   {\it Let $D$ be a 2-strong digraph of order $n\geq 3$ satisfying condition $(M)$. Suppose that
$\{u,v\}$ is a pair of non-adjacent vertices in $V(D)$ such that $d(u)+d(v)\leq 2n-2$. Then $D$ is Hamiltonian or $D$ contains a cycle of length $n-1$ passing through $u$ and avoiding $v$ (passing through $v$ and avoiding $u$).}\\

As an immediate corollary of   Theorems 3.6, 3.3 and Lemma 3.1, we obtain\\

\noindent\textbf{Corollary 3.7.} {\it Let $D$ be a 2-strong non-Hamiltonian digraph of order $n\geq 3$ satisfying condition $(M)$. Suppose that
$\{u,v\}$ is a pair of non-adjacent vertices in $V(D)$ such that $d(u)+d(v)\leq 2n-2$. Then $d(u)\leq n-1$,  $d(v)\leq n-1$ and $D$ contains at most one cycle of length two  passing through $u$  ($v$) }.

 \section{Preliminaries }

\noindent\textbf{Lemma 4.1.}  {\it  Let $D$ be a 2-strong  digraph of order $n\geq 3$ satisfying condition $(M)$. Suppose that $\{y,z\}$ is a pair of non-adjacent vertices in  $V(D)$ such that $d(y)+d(z)\leq 2n-2$ and
 $C=x_1x_2\ldots x_{n-k}x_1$ is a cycle in $D$ passing through $y$ and avoiding $z$, where $2\leq n-k\leq n-2$. If the subdigraph 
 $D\langle V(D)\setminus V(C)\rangle$ contains  a cycle
   passing through $z$ and $d(y,V(D)\setminus V(C))=0$,  
   then $D$ is Hamiltonian}.
 
\begin {proof}
Suppose, on the contrary, that 
 $D\langle V(D)\setminus V(C)\rangle$ contains  a cycle passing through $z$, but $D$ is not Hamiltonian. Since $D$ contains at most one cycle of length two passing through $y$ (Corollary 3.7), from $d(y,V(D)\setminus V(C))=0$ it follows that $d(y)\leq n-k$. Let $y_1y_2\ldots y_sy_1$ be a cycle through $z$ in $D\langle V(D)\setminus V(C)\rangle$, where $s\in [2,k]$. 
  
    By Theorem 3.3 we have that $D$ contains no cycle through $y$ and $z$. Therefore, for each pair of integers $i$ and $j$, 
where $i\in [1,n-k]$ and $j\in [1,s]$, 
$\overrightarrow{a} [x_i,y_j]+\overrightarrow{a}[y_{j-1},x_{i+1}]\leq 1$  (here, $y_0=y_s$ and $x_{n-k+1}=x_1$). This implies that for every $j\in [1,s]$ we have
$$
d^-(y_j, V(C))+d^+(y_{j-1}, V(C))=\sum_{i=1}^{n-k}(\overrightarrow{a} [x_i,y_j]+\overrightarrow{a}[y_{j-1},x_{i+1}]))
  \leq n-k. $$
Hence,
$$
d(y_1, V(C))+\cdots +d(y_s, V(C))=\sum_{j=1}^{s}(d^-(y_j, V(C))+d^+(y_{j-1}, V(C)))
  \leq s(n-k). \eqno (1)
$$
Since there is at most one cycle of length two through $z$ ($y$) (Corollary 3.7), it follows that for $A:= V(D)\setminus V(C)$ and for every $y_j\in \{y_1,\ldots , y_s\}\setminus \{z,y_1\}$ (we may assume that $y_1\not= z$) the following holds:
$$
d(z,A)\leq k, \quad d(y_1,A)\leq 2k-2 \quad \hbox{and} \quad d(y_j,A)\leq 2(k-2)+1=2k-3.
$$
Therefore, 
$$
d(y_1,A)+\cdots +d(y_s,A)\leq (s-2)(2k-3)+k+2k-2=2ks-3s-k+4. 
$$
Combining  this  with (1), we obtain
$$
d(y_1)+\cdots +d(y_s)\leq ns+ks-3s-k+4. 
$$
The last inequality together with  $d(y)\leq n-k$ implies that
 $$
d(y_1)+\cdots +d(y_s)+sd(y)\leq 2ns-3s-k+4. \eqno (2)
$$
 Notice that $\{y,y_1\}$, $\ldots$ , 
$\{y,y_s\}$ are $s$ distinct pairs of non-adjacent vertices. We will consider the cases when $s$ is even and $s$ is odd   separately.

Assume first that $s$ is even.  Using condition $(M)$ and (2), we obtain
$$
s(4n-3)/2\leq d(y_1)+\cdots +d(y_s)+sd(y)\leq 2ns-3s-k+4.
$$
Therefore, $2ns-1.5s\leq 2ns-3s-k+4$, i.e., $1.5s+k\leq 4$. The last inequality is impossible, since $k\geq s\geq 2$.

Assume next that $s$ is odd. Then $s\geq 3$. Since $d(y)\leq n-k$, and $d(z)\leq n-1$ by Corollary 3.7 (we may  assume that $z\not=y_s$), from condition $(M)$ it follows that $d(y)+d(y_s)\geq 2n+k-2$. Now, by condition $(M)$ and (2) we have, 
$$
(s-1)(4n-3)/2+2n+k-2\leq d(y_1)+\cdots +d(y_{s-1})+d(y_s)+sd(y)$$ $$\leq 2ns-3s-k+4.
$$ 
Hence, $$
2n(s-1 )-1.5(s-1)+2n+k-2\leq 2ns-3s-k+4.
$$ 
This means that $1.5s+2k\leq 4.5$, which is a contradiction. This contradiction completes the proof of Lemma 4.1. 
\end {proof}

\noindent\textbf{Lemma 4.2.} {\it Let $D$ be a 2-strong digraph of order $n\geq 3$ satisfying condition $(M)$. Suppose that $\{y,z\}$ is a pair of non-adjacent vertices in  $V(D)$ such that $d(y)+d(z)\leq 2n-2$ and
 $C=x_1x_2\ldots x_{n-2}zx_1$ is a cycle of length $n-1$  passing through $z$ and avoiding $y$ in $D$.  Then either $D$ is Hamiltonian or for every $k\in [2,n-3]$, the following holds:
$$
A(\{x_1,\ldots , x_{k-1}\}\rightarrow  \{x_{k+1},\ldots , x_{n-2}\})\not=\emptyset.
$$} 
\begin {proof}
Suppose that $D$ is not  Hamiltonian. Since $D$ is 2-strong, $n\geq 5$. Then  by Theorem 3.3, there is no cycle  through $y$ and $z$.
 Therefore, we have that if $x_i\rightarrow y$ with $i\in [1,n-3]$, then 
 $d^+(y,\{x_{i+1}, \ldots , x_{n-2}\})=0$
  (for otherwise, $x_1\ldots x_iyx_j\ldots x_{n-2}zx_1$, where $j\in [i+1,n-2]$, is a cycle through $y$ and $z$, a contradiction).
   Let $x_r\rightarrow y\rightarrow x_p$, $1\leq p<r\leq n-2$, and $p$, $r$ be chosen so that $p$ is minimal and $r$ is maximal with these properties.  Then
   $$
   d(y,\{x_{1}, \ldots , x_{p-1}\})=d(y,\{x_{r+1}, \ldots , x_{n-2}\})=0. \eqno (3)
   $$
   If $p=1$ and $r=n-2$, then by a similar argument as above, we conclude that if  $x_i\rightarrow z$ with $i\in [1,n-3]$, then 
$d^+(z,\{x_{i+1},\ldots , x_{n-2}\})=0$. Assume that $p\geq 2$ or $r\leq n-3$. Observe that $Q:=yx_p\ldots x_ry$ is a cycle  through $y$ which does not contain $z$, and $d(y,V(D)\setminus V(Q))=0$ because of (3). Therefore by Lemma 4.1, the subdigraph $D\langle V(D)\setminus V(Q)\rangle$ contains no cycle  through $z$ since $D$ is not Hamiltonian. This implies that
$$
d^-(z,\{x_1,\ldots , x_{p-1}\})=d^+(z,\{x_{r+1},\ldots , x_{n-2}\})=0
$$
since $x_{n-2}\rightarrow z\rightarrow x_1$. From the last equalities it follows that if there are $i, j$ such that  $x_i\rightarrow z$ and  $z\rightarrow x_j$ with $i<j$, then $i\geq p$, $j\leq r$ and $yx_p\ldots x_izx_j\ldots x_ry$ is a cycle passing through $y$ and  $z$, a contradiction. Thus, we may assume that for every pair of integers $i$ and $j$, $1\leq i< j\leq n-2$,
$$
\hbox{if} \quad x_i\rightarrow y, \quad \hbox{then} \quad yx_j\notin A(D)\quad \hbox{ and if} \quad x_i\rightarrow z, \quad \hbox{then} \quad zx_j\notin A(D).  \eqno (4)
$$
Now suppose that the theorem is not true. Then $D$ is not Hamiltonian and there  is an integer $k\in [2,n-3]$ such that
$$
A(\{x_1,\ldots , x_{k-1}\}\rightarrow  \{x_{k+1},\ldots , x_{n-2}\})=\emptyset.   \eqno (5)
$$
It is easy to see that there are  vertices $x_m$ and $x_l$ such that 
$y\rightarrow  x_m$, $z\rightarrow  x_l$ and 
$$
d^+(y, \{x_{m+1}, \ldots , x_{n-2}\})=d^+(z, \{x_{l+1}, \ldots , x_{n-2}\})=0. \eqno (6)
$$
Then by (4),
$$
d^-(y,\{x_1,\ldots , x_{m-1}\})=d^-(z,\{x_{1},\ldots , x_{l-1}\})=0. \eqno (7)
$$
Assume first that $m\leq l$. Since $D$ is 2-strong,  (4)
 and (7) imply  that $2\leq m\leq l\leq n-3$. Now from (5), (6) and (7)  it follows that:

(i) if $k\leq m$ or $k\geq l$, then (respectively)
$$
A(\{x_1,x_2,\ldots , x_{k-1}\}\rightarrow  \{y,z,x_{k+1},x_{k+2},\ldots , x_{n-2}\})=\emptyset  
$$
or
$$
A(\{y,z,x_1,x_2,\ldots , x_{k-1}\}\rightarrow  \{x_{k+1},x_{k+2},\ldots , x_{n-2}\})=\emptyset,   
$$
(ii) if $m < k< l$, then 
$
A(\{y,x_1,x_2,\ldots , x_{k-1}\}\rightarrow  \{z,x_{k+1},x_{k+2},\ldots , x_{n-2}\})=\emptyset   
$.
Thus, in each case we have that $D-x_k$ is not strong, which contradicts the condition that $D$ is 2-strongly connected.

Assume next that $m>l$. This case is similar to the first case and we omit the details. Lemma 4.2 is proved. \end {proof}  

The following lemma is proved by  \cite{[13]}. We present its proof for completeness.\\

\noindent\textbf{Lemma 4.3.} {\it Let $D$ be a 2-strong  digraph of order $n\geq 3$ satisfying condition $(M)$. Suppose that $\{y,z\}$ is a pair of non-adjacent vertices in  $V(D)$ such that $d(y)+d(z)\leq 2n-2$ and
 $C=x_1x_2\ldots x_{n-2}zx_1$ is a  cycle of length $n-1$ passing through $z$ and avoiding $y$ in $D$.  
If  $x_a\rightarrow x_b$ and there are  integers $l, s, f, t$ such that $1\leq l\leq a<s\leq f<b\leq t\leq n-2$ and  $\{x_f,x_t\}\rightarrow y\rightarrow \{x_l,x_s\}$, 
 then $D$ is Hamiltonian.}
 
 \begin {proof}
 Suppose, on the contrary, that $D$ is not Hamiltonian.  By Theorem 3.3, $D$ contains no cycle  through $y$ and $z$. Therefore, there are no integers $i$ and $j$, $1\leq i<j\leq n-2$, such that $x_i\rightarrow y \rightarrow x_j$ (for otherwise, $x_1\ldots x_iyx_j\ldots x_{n-2}zx_1$ is a cycle  through $y$ and $z$). Since the arcs $yx_l$, $yx_s$, $x_fy$, $x_ty$ are in $D$ and 
 $l\leq a< s\leq f<b\leq t$,
 it is easy to check that:
 
(i) if $z\rightarrow x_i$ with $i\in [a+1,f]$, then $C(y,z)=yx_l\ldots x_ax_b\ldots x_{n-2}zx_i\ldots x_fy$;

(ii) if $x_j\rightarrow z$ with $j\in [s,b-1]$, then $C(y,z)=x_1\ldots x_ax_b\ldots x_{t}yx_s\ldots x_jzx_1$.
Thus, in both cases we have a contradiction. Therefore,
 $$
 d^+(z,\{x_{a+1},\ldots , x_{f}\})= d^-(z,\{x_{s},\ldots , x_{b-1}\})=0,
$$
in particular, $d(z,\{x_{s},\ldots , x_{f}\})=0$ and the vertices $z$ and $x_s$ ($z$ and $x_f$) are not adjacent. The last equality together with the fact that  $D$ contains at most one cycle of length two passing through $z$ (Corollary 3.7) implies that
$$
d(z)=d(z,\{x_{1},\ldots , x_{s-1}\})+d(z,\{x_{f+1},\ldots , x_{n-2}\})\leq n+s-f-2. \eqno (8)
$$

Now we consider the  vertex $x_s$. It is not difficult to check that:

(iii) if $x_i\rightarrow x_s$ with $i\in [1,l-1]$, then
$C(y,z)=x_1\ldots x_ix_s\ldots x_{f}yx_l\ldots$ $ x_ax_b\ldots x_{n-2}zx_1$;

(iv) if $x_s\rightarrow x_j$ with
 $j\in [t+1,n-2]$,  then
$C(y,z)=x_1\ldots x_ax_b\ldots x_{t}yx_s$ $ x_j\ldots x_{n-2}zx_1$. In both cases we have a contradiction.
Therefore, we may assume that
 $$
d^-(x_s,\{x_1,\ldots , x_{l-1}\})=d^+(x_s,\{x_{t+1},\ldots , x_{n-2}\})=0.
$$
This implies that
$$
d(x_s)= d^+(x_s,\{x_{1},\ldots , x_{l-1}\})+ d^-(x_s,\{x_{t+1},\ldots , x_{n-2}\})
+ d(x_s,\{x_l, \ldots , x_t\})+
d(x_s,\{y\})$$ $$\leq  l-1+n-2-t+2(t-l+1)=n+t-l-1. \eqno (9)
$$
Without loss of generality, we may assume that $l$, $f$ are chosen as maximal as possible and $s$, $t$  are chosen as minimal as possible, i.e.,
$$
  d(y,\{x_{l+1},\ldots , x_{s-1}\})= d(y,\{x_{f+1},\ldots , x_{t-1}\})= 0.
$$
This, since $D$ contains at most one cycle of length two passing through $y$, implies that
$$
d(y)=d(y,\{x_{1},\ldots , x_{l}\})+ d(y,\{x_{s},\ldots , x_{f}\})+d(y,\{x_t,\ldots , x_{n-2}\})$$
$$\leq
l+f-s+1+n-2-t+2=n+l+f-s-t+1. 
$$
Since $\{y,z\}$ and $\{x_s,z\}$ are two distinct pairs of non-adjacent vertices, from (8), (9), the last inequality and condition $(M)$ it follows that
$$
4n-3\leq d(y)+2d(z)+d(x_s)\leq n+l+f-s-t+1+2n+2s-2f-4+n+t-l-1$$ $$=4n-4-(f-s)\leq 4n-4,
$$
which is a contradiction.
Lemma 4.3 is proved. 
\end {proof}

 \section{Proof of Theorem 1.12}
 
Recall the statement of Theorem 1.12.\\

\noindent\textbf{Theorem 1.12.} {\it Let $D$ be a 2-strong digraph of order $n\geq 3$ satisfying condition $(M)$. Then $D$ is Hamiltonian}.

\begin {proof}
By Theorem 3.4, the theorem is true if $D$ contains at most one pair of non-adjacent vertices. We may therefore assume that $D$ contains at least two distinct pairs of non-adjacent vertices. If the degrees sum of any two non-adjacent vertices at least $2n-1$, then  by  Meyniel's theorem, the theorem is true. We may therefore assume that $D$ contains a pair of non-adjacent vertices, say $y,z$, such that $d(y)+d(z)\leq 2n-2$. By Theorem 3.3, to prove the theorem, it suffices to prove that $D$ contains a cycle 
 through $y$ and $z$. If $d(y)+d(z)\geq 2n-3$, then by Lemma 3.5 we have that $D$ contains a cycle  trough $y$ and $z$, which, in turn, implies that $D$ is Hamiltonian (by Theorem 3.3). Thus, we may assume that $d(y)+d(z)\leq 2n-4$. By Theorem 3.6 we have that either $D$ is Hamiltonian or $D$ contains a cycle of length $n-1$ passing through $z$ and avoiding $y$ 
(passing through $y$ and avoiding $z$).

Suppose that $D$ is not Hamiltonian, i.e., $D$ contains no cycle  through $y$ and $z$. Let $C:=x_1x_2\ldots x_{n-2} zx_1$ be a cycle of length $n-1$ in $D$,  which does not contain $y$. Let $q$ be the maximum integer such that $y\rightarrow x_q$ and $k$ be the minimum integer such that $x_k\rightarrow y$.  Since $D$ is 2-strong and contains no cycle passing through $y$ and $z$, it follows that $k\geq q$ and there
  are some integers $p$,  $r$, $1\leq p< q\leq k<r\leq n-2$,  such that $x_r\rightarrow y\rightarrow x_p$  and  
 $$
 d(y,\{x_1,\ldots , x_{p-1}\})=d(y, \{x_{q+1}, \ldots , x_{k-1}\})= d(y, \{x_{r+1}, \ldots , x_{n-2}\})   
 $$
 $$
= d^-(y,\{x_p,\ldots , x_{q-1}\})=d^+(y,\{x_{k+1},\ldots , x_{r}\})=0. \eqno(10)
 $$
 Note that if $D$ contains a cycle of length two passing trough $y$, then $k=q$, otherwise $k>q$, $yx_k\notin A(D)$ and  $x_qy\notin A(D)$.  Therefore, it is not difficult to see that
 $$
d(y)= d^+(y,\{x_p,\ldots , x_{q}\})+d^-(y,\{x_{k},\ldots , x_{r}\})\leq q-p+r-k+2. \eqno(11)
 $$
 
 In order to prove the theorem, it is convenient for the digraph $D$ and the cycle $C$ to prove the following claims.

\noindent\textbf{Claim 5.1.} {\it If $p\geq 2$, then
 $d^-(x_{n-2},\{z,x_1,\ldots , x_{p-1}\})=0$}.
 
\begin {proof}
 Notice that $Q:=yx_p\ldots x_ry$ is a cycle passing through $y$ and avoiding $z$.  By (10) we have that $d(y,
V(D)\setminus V(Q))=0$. Now by Lemma 4.1 , the induced subdigraph 
 $D\langle V(D)\setminus V(Q)\rangle$ contains no cycle  through $z$. Then, since $x_{n-2}\rightarrow z\rightarrow x_1$, we have
$$
 d^-(z,\{x_1,\ldots , x_{p-1}\})=0 \quad \hbox{and} 
 \quad A(\{z,x_1,\ldots , x_{p-1}\}\rightarrow \{x_{r+1},\ldots , x_{n-2}\})=\emptyset.
 $$
The first equality together with 2-connectedness of $D$ implies that there is an integer $t\in [p,n-3]$ such that $x_t\rightarrow z$. The last equality means that  if $r\leq n-3$, then $d^-(x_{n-2},\{z,x_1,\ldots , x_{p-1}\})=0$. Assume that $r=n-2$, i.e., $x_{n-2}\rightarrow y$. In this case, we have that if $x_i\rightarrow x_{n-2}$ with $i\in [1,p-1]$ (respectively, $z\rightarrow x_{n-2}$), then $C(y,z)=x_1\ldots x_ix_{n-2}yx_p\ldots x_tzx_1$ (respectively, $C(y,z)=yx_p\ldots x_tzx_{n-2}y$), which is a contradiction. This proves that 
$d^-(x_{n-2},\{z,x_1,\ldots , x_{p-1}\})=0$. \end{proof} 

\noindent\textbf{Claim  5.2.} {\it Suppose that $k\geq q+1$ and 
$x_h \rightarrow x_l$, where $h\in [q,k-1]$ and $l\in [k+1,n-2]$. Then
$d^-(x_k,\{x_1,\ldots , x_{q-1}\})=0$}.

\begin {proof}
Assume that  Claim 5.2 is not true. Then for some $i\in [1,q-1]$, $x_i \rightarrow x_k$. Then, since the arcs $yx_q$, $x_ky$, $x_hx_l$ are in $D$ and $i< q\leq h<k<l$, we have a cycle $C(y,z)=x_1\ldots x_ix_kyx_q\ldots x_hx_l\ldots x_{n-2}zx_1$, which contradicts our initial supposition. \end{proof} 

\noindent\textbf{Claim 5.3.} {\it Suppose that $k\geq q+1$, 
$x_h\rightarrow x_l$ with $h\in [q,k-1]$ and $l\in [k+1,r]$ (possibly, $r=n-2$). Then there is an integer $f\geq 0$ such that $l+f\leq r$, 
$x_{l+f}\rightarrow y$, $d(y,\{x_l,\ldots ,x_{l+f-1}\})=0$ 
(possibly, 
$\{x_l,\ldots ,x_{l+f-1}\}=\emptyset$). Moreover,  either there is a vertex $x_g$ with $g\in [l+f+1,n-2]$ such that $x_k\rightarrow x_g$ or  there is a vertex $x_{c}$ with $c\in [k,l-1]$ such that
 $x_{c}\rightarrow z$}.

\begin {proof}
By Claim 5.2,
$$d^-(x_k,\{x_1,\ldots , x_{q-1}\})=0. \eqno (12) $$
Since $l\leq r$ and  $x_{r}\rightarrow y$, obviously there is an integer 
$f\geq 0$ such that $l+f\leq r$, 
$x_{l+f}\rightarrow y$, $d^-(y,\{x_l,\ldots ,x_{l+f-1}\})=0$ 
 $($possibly $\{x_l,\ldots , x_{l+f-1}\}=\emptyset$). This together with \\ $d^+(y,\{x_l,\ldots ,x_{l+f-1}\})=0$ implies that 
$$d(y,\{x_l,\ldots ,x_{l+f-1}\})=0. \eqno (13)$$
Now suppose that the claim is not true. Then 
 $$
 d^+(x_k,\{ x_{l+f+1}, \ldots , x_{n-2}\})=0 \\\ \hbox{and 
 }  \\\ d^-(z,\{x_k,\ldots , x_{l-1}\})=0. \eqno (14)
 $$   
The   second equality of (14) together with $d^+(y,\{x_k,\ldots ,x_{l-1}\})=0$ and the fact that there is no path of length two between $y$ and $z$ (Lemma 3.5) implies that the vertices $x_k$, $z$ are not adjacent and
$$
d(z,\{x_k,\ldots , x_{l-1}\})+d(y,\{x_k,\ldots , x_{l-1}\})\leq l-k.
$$
This together with (13), (10) and  the fact that there is at most one cycle of length two  through $z$ (Corollary 3.7) implies that 
$$
d(y)+d(z)=d^+(y,\{x_p,\ldots , x_{q}\})+d(y,\{x_k,\ldots , x_{l-1}\})+
d(z,\{x_k,\ldots , x_{l-1}\})$$
 $$
 +d^-(y,\{x_{l+f},\ldots , x_{r}\})+
d(z,\{x_1,\ldots , x_{k-1}\})+d(z,\{x_l,\ldots , x_{n-2}\})$$
 $$
 \leq q-p+1+l-k+r-l-f+1+k-1+n-2-l+2$$ $$= n+q+r+1-p-l-f. 
$$
Now consider the vertex $x_k$. Note that $d(x_k,\{y\})=1$ since $k\geq q+1$. Using (12) and the first equality of  (14), we obtain
$$
d(x_k)=d^+(x_k,\{x_1,\ldots , x_{q-1}\})+d(x_k,\{x_q,\ldots , x_{l+f}\})+
d^-(x_k,\{x_{l+f+1},\ldots , x_{n-2}\})$$ $$+d^+(x_k,\{y\})\leq q-1+2l+2f-2q+n-2-l-f+1=n+l+f-q-2.
 $$
Combining the last two inequalities,  $d(z)\leq n-1$ (Corollary 3.7) and $r\leq n-2$, we obtain 
$$
d(y)+d(z)+d(x_k)+d(z)\leq 3n+r-p-2\leq 4n-4-p,
$$
which contradicts condition $(M)$, since
$\{y,z\}$, $\{z,x_k\}$ are two distinct pairs of non-adjacent vertices. This contradiction completes the proof of  Claim 5.3.
\end{proof} 

\noindent\textbf{Claim 5.4.} {\it If $p\geq 2$, then 
$A(\{x_1,\ldots , x_{p-1}\}\rightarrow \{x_{k+1},\ldots , x_{n-2}\})=\emptyset.
$}

\begin {proof}
 Suppose, on the contrary, that $p\geq 2$ and $x_a\rightarrow x_b$ with $a\in [1,p-1]$ and $b\in [k+1,n-2]$. 
 Let $b$ be the maximum with these properties, i.e.,
$$
A(\{x_1,\ldots , x_{p-1}\}\rightarrow \{x_{b+1},\ldots , x_{n-2}\})=\emptyset. \eqno (15)
$$
Notice that $Q:=yx_p\ldots x_ry$
is a cycle in $D$ and $d(y,V(D)\setminus V(Q))=0$ by (10). Therefore by Lemma 4.1, the subdigraph $D\langle V(D)\setminus V(Q)\rangle$ does not contain a cycle through $z$. In particular, 
$$
d^-(z,\{x_1,\ldots , x_{p-1}\})=0, \eqno (16)
$$
 and if $r\leq n-3$, then 
$$
d^+(z,\{x_{r+1},\ldots , x_{n-2}\})=0 \quad \hbox{and} \quad 
A(\{x_1,\ldots , x_{p-1}\}\rightarrow \{x_{r+1},\ldots , x_{n-2}\})=\emptyset.
 \eqno (17)
 $$
By Claim 5.1, we have 
$$
d^-(x_{n-2},\{z,x_{1},\ldots , x_{p-1}\})=0. \eqno (18)
$$
From (17) and (18) it follows that $b\leq r$ and, if $r=n-2$, then
 $b\leq n-3$. In both cases we have that $b\leq n-3$.
 
  If $x_i\rightarrow z$ with $i\in [p,b-1]$, then
 $C(y,z)=x_1\ldots x_ax_b\ldots x_ryx_p\ldots x_izx_1$, a contradiction. We may therefore assume that 
 $d^-(z,\{x_{p},\ldots , x_{b-1}\})=0$. This together with (16) implies that 
$$
d^-(z,\{x_{1},\ldots , x_{b-1}\})=0. \eqno (19)
$$
Applying Lemma 4.2 to the vertex $x_b$, we obtain that 
$$
A(\{x_1,\ldots , x_{b-1}\}\rightarrow \{x_{b+1},\ldots , x_{n-2}\})\not=\emptyset.
$$
 Let $x_s\rightarrow x_t$, where $s\in [1,b-1]$ and $t\in [b+1,n-2]$. Choose  $t$  maximal with these properties, i.e.,
$$
A(\{x_1,\ldots , x_{b-1}\}\rightarrow \{x_{t+1},\ldots , x_{n-2}\})=\emptyset. \eqno (20)
$$
From (15) it follows that $s\geq p$, i.e., $s\in [p,b-1]$. If $x_i\rightarrow y$ with $i\in [b,t-1]$, then
 $C(y,z)=x_1\ldots x_ax_b\ldots x_iyx_p\ldots x_sx_t\ldots x_{n-2}zx_1$, a contradiction. We may therefore assume that \\
 $d^-(y,\{x_{b},\ldots , x_{t-1}\})=0$. This  together with 
 $d^+(y,\{x_b,\ldots , x_{t-1}\})=0$ implies that 
 $$
 d(y,\{x_b,\ldots , x_{t-1}\})=0. \eqno (21)
 $$
 In particular, the vertices $x_b$ and $y$ are not adjacent,  $t\leq r$ and $b\leq r-1$ since $b+1\leq t\leq r$  (i.e., $A(\{x_p,\ldots , x_{b-1}\}\rightarrow \{x_{r+1},\ldots , x_{n-2}\})=\emptyset$). 
 Using Lemma 4.3, we obtain 
 $$
A(\{x_p,\ldots , x_{q-1}\}\rightarrow \{x_{k+1},\ldots , x_{r}\})=\emptyset \quad \hbox{and} \quad d^-(x_{k+1},\{x_p,\ldots x_{q-1}\})=0. \eqno (22)
$$
Then, since $t\leq r$ and (20), we have that
 $
A(\{x_p,\ldots , x_{q-1}\}\rightarrow \{x_{b+1},\ldots , x_{n-2}\})=\emptyset$.  
This together with (15) implies that
$$
A(\{x_1,\ldots , x_{q-1}\}\rightarrow \{x_{b+1},\ldots , x_{n-2}\})=\emptyset. 
$$
Therefore, $s\geq q$, i.e., $s\in [q,b-1]$. Since $b\leq r-1$, and  $x_b$, $y$ are not adjacent, there is an integer $f\geq 0$ such that  
$d^-(y,\{x_b,\ldots , x_{b+f}\})=0$ and $x_{b+f+1}\rightarrow y$. Then, since (21) and   $d^+(y,\{x_b,\ldots , x_{b+f}\})=0$ we have that $t\leq b+f+1$ and 
$$
d(y,\{x_b,\ldots , x_{b+f}\})=0. \eqno (23)
$$
This together with (10) implies that 
$$
d(y)=d^+(y,\{x_p,\ldots , x_{q}\})+d^-(y,\{x_k,\ldots , x_{b-1}\})+d^-(y,\{x_{b+f+1},\ldots , x_{r}\})$$ $$\leq q-p+1+b-k+r-b-f=q+r+1-p-k-f.  \eqno   (24)
$$
From (19), $d^+(y,\{x_k,\ldots , x_{b-1}\})\leq 1$ and the fact that there is no path of length two between $y$ and $z$ (Lemma 3.5) it follows that 
$$
d(y,\{x_k,\ldots , x_{b-1}\})+d(z,\{x_k,\ldots , x_{b-1}\})\leq b-k+1.  
$$
This together with (10), (23) and the fact that there is at most one cycle of length two  through $z$ (Corollary 3.7) implies that
$$
d(y)+d(z)=d^+(y,\{x_p,\ldots , x_{q}\})+d(y,\{x_k,\ldots , x_{b-1}\})+
d(z,\{x_k,\ldots , x_{b-1}\})$$ $$+d^-(y,\{x_{b+f+1},\ldots , x_{r}\})+
d(z,\{x_1,\ldots , x_{k-1}\})+d(z,\{x_b,\ldots , x_{n-2}\})$$  
$$\leq q-p+1+b-k+1+r-b-f+k-1+n-2-b+2$$ $$= n+1+q-p+r-b-f.  \eqno (25)
$$
 Since $t\leq b+f+1$ and (20), it follows that 
 $$
A(\{x_p,\ldots , x_{b-1}\}\rightarrow \{x_{b+f+2},\ldots , x_{n-2}\})=\emptyset. \eqno (26)
$$
In particular, from $b\geq k+1$ and (26) it follows that
$$
d^+(x_k,\{x_{b+f+2},\ldots , x_{n-2}\})=0. \eqno (27)
$$
We will consider the cases  $b\geq k+2$, $b=k+1$ separately.\\

\noindent\textbf{Case 1.} $b\geq k+2$.

Then by the first equality of (22) we have 
$$
d^-(x_{b-1},\{x_{p},\ldots , x_{q-1}\})=0. \eqno (28)
$$
Using the fact that there is no path of length two between $y$ and $z$ (Lemma 3.5) and (19), we obtain that $d(x_{b-1},\{y,z\})\leq 1$. This together with $d^+(x_{b-1},\{x_{b+f+2},\ldots , x_{n-2}\})=0$  (by (26))
 and (28) implies that 
$$
d(x_{b-1})=d(x_{b-1},\{x_{1},\ldots , x_{p-1}\})+d^+(x_{b-1},\{x_{p},\ldots , x_{q-1}\})+d(x_{b-1},\{x_{q},\ldots , x_{b+f+1}\})$$ $$ +d^-(x_{b-1},\{x_{b+f+2},\ldots , x_{n-2}\})+d(x_{b-1},\{y,z\})
  \leq  2p-2+q-p$$ $$+2b+2f+2-2q+n-2-b-f-1+1=n+p-q+b+f-2. \eqno (29)
$$
Now we divide this case into the following subcases.\\

\noindent\textbf{Subcase 1.1.} {\it The vertices $x_{b-1}$ and $y$ are not adjacent}.

Then $\{y,x_{b-1}\}$ and $\{y,z\}$ are two distinct pairs of non-adjacent vertices. Since $p\geq 2$, $r\leq n-2$, $f\geq 0$ and $k\geq q$, combining (25), (24) and (29), we obtain 
$$
d(y)+d(z)+d(y)+d(x_{b-1})\leq n+1+q-p+r-b-f+q+r+1-p-k-f$$ $$+n+p-q+b+f-2
=2n+2r+q-p-k-f\leq 4n-4-(k-q)-f-p,
 $$
which contradicts condition $(M)$.\\

\noindent\textbf{Subcase 1.2.} {\it The vertices $x_{b-1}$ and $y$ are  adjacent}.
  
  Then $x_{b-1}\rightarrow  y$. Therefore by Lemma 3.5 and (19), the vertices $z$ and $x_{b-1}$ are not adjacent. Since $d(z)\leq n-2$ (because of $d(z,\{y,x_{b-1}\})=0$ and Corollary 3.7) and $r\leq n-2$, from (25) and (29) it follows that 
  $$
  d(y)+d(z)+d(x_{b-1})+d(z)\leq n+1+q-p+r-b-f+n+p-q+b+f-2+n-2$$ $$=3n-3 +r\leq 4n-5,
  $$
  which contradicts condition $(M)$. The discussion of Case 1 is completed. \\
  
\noindent\textbf{Case 2.} $b=k+1$. 

We divide this case into the following subcases.\\

\noindent\textbf{Subcase 2.1.} $s\leq k-1$.

Then $k\geq q+1$ since $s\geq q$. Then $yx_k\notin A(D)$ by the definition of $q$ and $k$. Recall that the vertices   $z$, $x_k$
are not adjacent by (19) and Lemma 3.5. Now it is easy to see that $d(z)\leq n-2$. Since $x_{s}\rightarrow  x_t$ with $s\in [q,k-1]$ and $t\in [b+1,n-2]$, by Claim 5.2 we have that 
$d^-(x_k,\{x_1,\ldots , x_{q-1}\})=0$. This together with (27) and  $b=k+1$ implies that 
$$
d(x_k)=d^+(x_k,\{x_1,\ldots , x_{q-1}\})+d(x_k,\{x_q,\ldots , x_{b+f+1}\})+d^-(x_k,\{x_{b+f+2},\ldots , x_{n-2}\})$$ 
$$+d^+(x_k,\{y\})\leq q-1+2b+2f+2-2q+n-2-b-f-1+1=n+k-q+f.  
$$
This together with (24) and $d(z)\leq n-2$, we obtain 
$$
 d(y)+d(x_k)+2d(z)\leq q+r+1-p-k-f+n+k-q+f+2n-4$$ $$=3n+r-p-3\leq 4n-5-p,
$$
which is a contradiction since $\{y,z\}$ and $\{x_k,z\}$ are two distinct pairs of non-adjacent vertices. \\

\noindent\textbf{Subcase 2.2.} $s= k$.

From $b=k+1$, $t\in [b+1=k+2,b+f+1]$ and (23) it follows that 
$$
d(y,\{x_{k+1},\ldots , x_{t-1}\})=0, \eqno (30)
$$
in particular, the vertices $y$ and $x_{k+1}$ are not adjacent. Observe that $R:=yx_p\ldots x_kx_t \ldots \\ x_ry$ is a cycle in $D$ passing through $y$, avoiding $z$ and $d(y,V(D)\setminus V(R))=0$. By Lemma 4.1, the induced subdigraph $D\langle V(D)\setminus V(R)\rangle$ contains no cycle  through $z$. In particular, this means that
$$
A(\{x_{k+1},\ldots , x_{t-1}\}\rightarrow \{x_{r+1},\ldots , x_{n-2}\})=\emptyset, \\\ \hbox{hence} \\\ d^+(x_{k+1},\{x_{r+1},\ldots , x_{n-2}\})=0, \eqno (31)
$$
for otherwise, if $x_i\rightarrow x_j$ with $i\in [k+1,t-1]$ and $j\in [r+1,n-2]$, then $H:=x_1\ldots x_ax_{k+1}\ldots x_ix_j \\ \ldots  x_{n-2}zx_1$ is a cycle in $D\langle V(D)\setminus V(R)\rangle$ through $z$, a contradiction.\\

\noindent\textbf{Subcase 2.2.1.} There is an integer  $l\in [b+f+2,n-2]$ such that $x_{k+1}\rightarrow x_l$ and
$$
d^+(x_{k+1},\{x_{l+1},\ldots , x_{n-2}\})=0. \eqno (32)
$$
Then $b+f+2\leq n-2$, and $l\leq r$ because of the first equality of (31). Recall that $t\leq b+f+1\leq l-1$. Hence, $l\geq t+1$. If $x_i\rightarrow z$ with $i\in [t,l-1]$, then $C(y,z)=x_1\ldots x_ax_{k+1}x_l\ldots x_ryx_q\ldots x_kx_t\ldots x_izx_1$, a contradiction. We may therefore assume that 
 $d^-(z,\{x_{t},\ldots , x_{l-1}\})=0$. This together with \\
$d^+(y,\{x_{t},\ldots , x_{l-1}\})=0$ and the fact that there is no path of length two between $y$ and $z$ implies that 
$$
d(y,\{x_{t},\ldots , x_{l-1}\})+d(z,\{x_{t},\ldots , x_{l-1}\})\leq l-t. 
$$
Combining this, (10) and (30), we obtain   
$$
d(y)+d(z)=
d^+(y,\{x_{p},\ldots , x_{q}\})+d^-(y,\{x_{k}\})+ 
d(y,\{x_{t},\ldots , x_{l-1}\})+d(z,\{x_{t},\ldots , x_{l-1}\})$$ $$+
d^-(y,\{x_{l},\ldots , x_{r}\})+d(z,\{x_{1},\ldots , x_{t-1}\})+
d(z,\{x_{l},\ldots , x_{n-2}\})$$ $$\leq q-p+1+1+l-t+r-l+1+t-1+n-2-l+2$$ $$\leq n+2+q+r-p-l.  \eqno (33)
$$
For the vertex $x_{k+1}$, using  (32) and the second equality of (22), we obtain 
$$
d(x_{k+1})=
d(x_{k+1},\{x_{1},\ldots , x_{p-1}\})+d^+(x_{k+1},\{x_{p},\ldots , x_{q-1}\})+ 
d(x_{k+1},\{x_{q},\ldots , x_{l}\})$$ $$+d^-(x_{k+1},\{x_{l+1},\ldots , x_{n-2}\})+
d(x_{k+1},\{z\})$$ $$\leq 2p-2+q-p+2l-2q+n-2-l+2=n-2+p-q+l. 
$$
This together with (33), (24),  $r\leq n-2$, $k\geq q$ and $p\geq 2$ implies that 
$$
d(y)+d(z)+d(y)+ d(x_{k+1})\leq n+2+q+r-p-l+q+r+1-p-k-f+n-2+p-q+l$$ $$=2n+1+q+2r-p-k-f\leq 4n-3-(k-q)-p-f\leq 4n-5,
$$
which contradicts condition $(M)$  since $\{y,z\}$ and $\{y,x_{k+1}\}$ are two distinct pairs of non-adjacent vertices.\\

\noindent\textbf{Subcase 2.2.2.} There is no $l\in [b+f+2,n-2]$ such that $x_{k+1}\rightarrow x_l$.

 Then 
 $d^+(x_{k+1},\{x_{b+f+2},\ldots , x_{n-2}\})=0$. This together with the second equality of (22) implies that
$$
d(x_{k+1})=
d(x_{k+1},\{x_{1},\ldots , x_{p-1}\})+d^+(x_{k+1},\{x_{p},\ldots , x_{q-1}\})$$ $$+ 
d(x_{k+1},\{x_{q},\ldots , x_{b+f+1}\})+d^-(x_{k+1},\{x_{b+f+2},\ldots , x_{n-2}\})+
d(x_{k+1},\{z\})$$ $$\leq 2p-2+q-p+2b+2f+2-2q+n-2-b-f-1+2$$ $$=n-1+p-q+b+f.  
$$
Combining this, $b=k+1$, (24) and $d(z)\leq n-2$, we obtain
$$ 
2d(y)+ d(x_{k+1})+d(z)\leq 2q+2r+2-2p-2k-2f+n-1+p-q+b+f$$ $$+ n-2= 2n+q+2r-p-k-f\leq 4n-4-(k-q)-p-f,
$$
which contradicts condition $(M)$. In each case we obtain a contradiction and hence the discussion of Case 2 is completed. This completes the proof of Claim 5.4. \end{proof} 

Now we are ready to complete the proof of the main result.

By Claim 5.4, if $p\geq 2$, then
$
A(\{x_1,\ldots , x_{p-1}\}\rightarrow \{x_{k+1},\ldots , x_{n-2}\})=\emptyset. 
$
Similarly,  if $r\leq n-3$, then 
$
A(\{x_1,\ldots , x_{q-1}\}\rightarrow \{x_{r+1},\ldots , x_{n-2}\})=\emptyset$. 
Using Lemma 4.3, we obtain
$
A(\{x_p,\ldots , x_{q-1}\}\rightarrow \{x_{k+1},\ldots , x_{r}\})=\emptyset. 
$
From the last three equalities  it follows that
$$
A(\{x_1,\ldots , x_{q-1}\}\rightarrow \{x_{k+1},\ldots , x_{n-2}\})=\emptyset. \eqno (34)
$$
From (34) and Lemma 4.2 it follows that $k\geq q+1$. Applying Lemma 4.2 to the vertices $x_q$ and $x_k$, we obtain
$$
A(\{x_1,\ldots , x_{q-1}\}\rightarrow \{x_{q+1},\ldots , x_{n-2}\})\not=\emptyset, \\\   
 A(\{x_1,\ldots , x_{k-1}\}\rightarrow \{x_{k+1},\ldots , x_{n-2}\})\not=\emptyset.$$ 
 Let $x_a\rightarrow x_b$ and
 $x_h\rightarrow x_l$ with $a\in [1,q-1]$,  $b\in [q+1,n-2]$, $h\in [1,k-1]$ and $l\in [k+1,n-2]$. Choose  $b$ maximal and $h$ minimal with these properties, i.e.,  
$$
A(\{x_1,\ldots , x_{q-1}\}\rightarrow \{x_{b+1},\ldots , x_{n-2}\})= 
A(\{x_1,\ldots , x_{h-1}\}\rightarrow \{x_{k+1},\ldots , x_{n-2}\})
=\emptyset. \eqno (35)
$$ 
From (34) it follows that $b\leq k$ and $h\geq q$, i.e., $b\in [q+1,k]$ and $h\in [q,k-1]$. 
  If $h\leq b-1$, then $C(y,z)=x_1\ldots 
  x_ax_b\ldots x_kyx_q\ldots x_hx_l\ldots x_{n-2}zx_1$, a contradiction. We may therefore assume that $h\geq b$, which in turn implies that $k\geq q+2$.
By Lemma 4.2, $A(\{x_1,\ldots , x_{b-1}\}\rightarrow \{x_{b+1},\ldots , x_{n-2}\})\not 
=\emptyset$. Let 
 $x_s\rightarrow x_t$, where $s\in [1,b-1]$ and $t\in [b+1,n-2]$. Choose $t$ maximal with this property, i.e.,
$$
A(\{x_1,\ldots , x_{b-1}\}\rightarrow \{x_{t+1},\ldots , x_{n-2}\})=\emptyset. \eqno (36)
$$   
From (35) it follows that $s\geq q$ and $t\leq k$, i.e., $s\in [q,b-1]$ and $t\in [b+1,k]$. We may assume that $l$ (recall that $x_h\rightarrow x_l$, $l\geq k+1$) is chosen so that   
  $$
  d^+(x_h,\{x_{k+1}, \ldots , x_{l-1}\})=0. \eqno (37)
  $$
  We consider the cases $l\leq r$ and $l\geq r+1$ separately.\\
  
 \noindent\textbf{Case 1.} $l\leq r$. 

For this case, it is not difficult to check that the conditions of Claim 5.3 hold.  Therefore, there is an integer $f\geq 0$ such that $l+f\leq r$, 
$x_{l+f}\rightarrow y$, $d(y,\{x_l,\ldots , x_{l+f-1}\})=0$ (possibly, $\{x_l,\ldots , x_{l+f-1}\}=\emptyset$), and either there is a vertex $x_g$ with $g\in [l+f+1,n-2]$ such that $x_k\rightarrow x_g$ or there is a vertex $x_{c}$ with $c\in [k,l-1]$ such that
 $x_{c}\rightarrow z$.

Assume first that $t\geq h+1$.
Then, since the arcs $yx_q$, $x_ax_b$, $x_sx_t$, $x_hx_l$, $x_ky$, $x_{l+f}y$  are in $D$ and $1\leq a\leq q-1<s<b\leq h<t\leq k<l\leq l+f\leq r\leq n-2$, we have that 
 $C(y,z)=x_1\ldots x_ax_b\ldots x_hx_l\ldots x_{l+f}yx_q\ldots x_sx_t\ldots x_{c}zx_1$, or 
 $C(y,z)=x_1\ldots x_ax_b\ldots x_hx_l\ldots x_{l+f}yx_q\ldots\\ x_sx_t\ldots  x_kx_g\dots x_{n-2}zx_1$ when $x_{c} \rightarrow z$ or when $x_k \rightarrow x_g$ respectively. In each case we have a contradiction.
 
 Assume next that $t\leq h$. By Lemma 4.2, $A(\{x_1,\ldots , x_{t-1}\}\rightarrow \{x_{t+1},\ldots , x_{n-2}\})\not 
=\emptyset$. Let 
 $x_{s_1}\rightarrow x_{t_1}$, where $s_1\in [1,t-1]$ and $t_1\in [t+1,n-2]$. Choose $t_1$ maximal with this property, i.e., 
$$
A(\{x_1,\ldots , x_{t-1}\}\rightarrow \{x_{t_1+1},\ldots , x_{n-2}\})=\emptyset. \eqno (38)
$$ 
From (36) (respectively, from (35)) it follows that $s_1\geq b$, i.e., $s_1\in [b,t-1]$ (respectively, $t_1\leq k$, i.e., $t_1\in [t+1,k]$). 
If $t_1\geq h+1$, then $C(y,z)=x_1\ldots x_ax_b\ldots x_{s_1}x_{t_1}\ldots x_kyx_q\ldots x_sx_t\ldots x_hx_l\dots x_{n-2}zx_1$,  a contradiction. 
We may therefore assume that $t_1\leq h$.   By Lemma 4.2, 
$$
A(\{x_1,\ldots , x_{t_1-1}\}\rightarrow \{x_{t_1+1},\ldots , x_{n-2}\})\not 
=\emptyset.
$$
 Let 
 $x_{s_2}\rightarrow x_{t_2}$, where $s_2\in [1,t_1-1]$ and $t_2\in [t_1+1,n-2]$. Choose $t_2$ maximal with this property, i.e.,
$$
A(\{x_1,\ldots , x_{t_1-1}\}\rightarrow \{x_{t_2+1},\ldots , x_{n-2}\})=\emptyset.  
$$ 
From (38) (respectively, from (35)) it follows that $s_2\geq t$, i.e., $s_2\in [t,t_1-1]$ (respectively, $t_2\leq k$, i.e., $t_2\in [t_1+1,k]$).
 
 Assume first that $t_2\geq h+1$. Then it is not difficult to see that $C(y,z)=x_1\ldots x_ax_b\ldots x_{s_1}x_{t_1}\ldots \\ x_hx_l\ldots x_{l+f}yx_q\ldots$ $ x_{s}x_{t}\ldots x_{s_2}x_{t_2}\ldots x_{c}zx_1$
or $C(y,z)=x_1\ldots x_ax_b\ldots x_{s_1}x_{t_1}\ldots x_hx_l\ldots x_{l+f}y \\ x_q\ldots  x_{s}x_{t}\ldots x_{s_2}x_{t_2}\ldots x_kx_g\dots$ $ x_{n-2}zx_1$ when $x_{c} \rightarrow z$ or when $x_k \rightarrow x_g$, respectively. In each case we have a contradiction.

 Continuing this process, we finally conclude that for some $m\geq 0$,
 $t_m\in [h+1,k]$ (here, $t_0=t$) since all the vertices $x_{t},x_{t_1},\ldots , x_{t_m}$ are distinct and  in $\{x_{q+1},\ldots , x_k$\}.  We already have constructed   a cycle $C(y,z)$ when $m\in \{0,1,2\}$. Assume that $m\geq 3$.  By the above arguments we have that:
 
  If $m\geq 3$ is odd, then
$C(y,z)=x_1\ldots x_ax_b\ldots x_{s_1}x_{t_1}\ldots x_{s_m}x_{t_m}\ldots x_kyx_q \ldots x_{s}x_{t}\ldots x_{s_2}x_{t_2}\ldots \\ $ $ x_{s_{m-1}}x_{t_{m-1}} \ldots x_{h}x_{l}\ldots  x_{n-2}zx_1$.

 If $m\geq 4$ is even, then 
  $C(y,z)=x_1\ldots x_{a}x_{b}\ldots x_{s_1}x_{t_1}\ldots x_{s_{m-1}}x_{t_{m-1}}\dots  x_{h}x_{l}\ldots x_{l+f}yx_q \ldots  x_{s}x_{t}\\\ldots  x_{s_2}x_{t_2} \ldots x_{s_{m}}x_{t_{m}}\ldots x_{c}zx_1$
or
$C(y,z)=x_1\ldots x_{a}x_{b}\ldots x_{s_1}x_{t_1}\ldots$ $  x_{s_{m-1}}x_{t_{m-1}}\ldots   x_{h}x_{l} \ldots\\$ $ x_{l+f}yx_q\ldots x_{s}x_{t}\ldots  x_{s_2}x_{t_2}\ldots x_{s_m}x_{t_m}\ldots x_kx_g\ldots x_{n-2}zx_1$ when $x_{c}\rightarrow z$   or when $x_k\rightarrow x_g$, respectively.  In all cases we have a cycle  through $y$ and $z$, which contradicts our supposition and hence the discussion of Case 1 is completed.\\
 
 \noindent\textbf{Case 2.} $l\geq r+1$.

Then $r\leq n-3$. 
Recall that $h\in [b,k-1]$, $x_h\rightarrow x_l$ and $x_s\rightarrow x_t$, where $l\leq n-2$, $s\in [q,b-1]$
and $t\in [b+1,k]$. Note that $\{y,x_h\}$, $\{y,z\}$ are two distinct pairs of non-adjacent vertices.\\

\noindent\textbf{Subcase 2.1.} $t\geq h+1$.

 Since $s\in [q,b-1]$
and $t\in [h+1,k]$, we have that $Q:=yx_p\ldots x_sx_t\ldots x_ry$ is a cycle in $D$ and $d(y,V(D)\setminus V(Q))=0$. If $a\leq p-1$, then 
$H:=x_1\ldots x_ax_b\ldots x_hx_l\ldots x_{n-2}zx_1$ is a cycle in
 $D\langle V(D)\setminus V(Q)\rangle$ passing through $z$, which contradicts Lemma 4.1.
  We may therefore assume
that $a\geq p$, i.e., $a\in [p,q-1]$.

Assume first that $b\leq h-1$. Then $q+1\leq b\leq h-1\leq k-2$ and $k\geq q+3$. From the first equality of (35) it follows that  
$d^-(x_h,\{x_1,\ldots , x_{q-1}\})=0$. This equality together with (37) implies that
$$
 d(x_h)=d^+(x_h,\{x_1,\ldots , x_{q-1}\})+d(x_h,\{x_q,\ldots , x_k\})+d^-(x_h,\{x_{k+1},\ldots , x_{l-1}\})$$ $$+d(x_h,\{x_{l},\ldots , x_{n-2}\})+d(x_h,\{z\})\leq q-1+2k-2q+l-1-k+2n-2l-2+2$$ $$
 = 2n-2-q+k-l. 
 $$
This together with (11) and $d(z)\leq n-1$ implies that
$$
 2d(y)+d(x_h)+d(z)\leq 2q-2p+2r-2k+4+2n-2-q+k-l+n-1$$ $$\leq 4n-2+(r-l)+(q-k)-2p,
 $$
which contradicts condition $(M)$.

Assume that $b=h$, i.e., $x_a\rightarrow x_h$. We may assume that $a$ is chosen so that 
$d^-(x_h,\{x_1,\ldots , x_{a-1}\})=0$. This and  (37) imply that
$$
 d(x_h)=d^+(x_h,\{x_1,\ldots , x_{a-1}\})+d(x_h,\{x_a,\ldots , x_k\})+d^-(x_h,\{x_{k+1},\ldots , x_{l-1}\})$$ $$+d(x_h,\{x_{l},\ldots , x_{n-2}\})+d(x_h,\{z\})\leq a-1+2k-2a+l-1-k+2n-2l-2+2$$ $$
 = 2n-2-a+k-l.  \eqno (39)
 $$
 Since $a\geq p$, it is not difficult to check that if $z\rightarrow x_i$ with $i\in [a+1,s]$, then $C(y,z)=yx_p\ldots x_ax_hx_l\ldots\\  x_{n-2}zx_i\ldots  x_sx_t\ldots x_ky$, which is a contradiction. Therefore,  $d^+(z,\{x_{a+1},\ldots , x_{s}\})=0$. This together with  $d^-(y,\{x_{a+1},\ldots , x_{s}\})=0$ and the fact that there is no path of length two between $y$ and $z$ implies that
$$
d(y,\{x_{a+1},\ldots , x_{s}\})+d(z,\{x_{a+1},\ldots , x_{s}\})\leq s-a.
$$ 
Using this and (10), we obtain 
$$
d(y)+d(z)= 
d^+(y,\{x_{p},\ldots , x_{a}\})+
d(y,\{x_{a+1},\ldots , x_{s}\})+d(z,\{x_{a+1},\ldots , x_{s}\})$$ $$+
d^-(y,\{x_{k},\ldots , x_{r}\})+d(z,\{x_{1},\ldots , x_{a}\})
+d(z,\{x_{s+1},\ldots , x_{n-2}\})$$ $$\leq a-p+1+s-a+r-k+1+a+n-2-s+1=n+1+a-p+r-k.
$$
Combining this, (11) and (39), we obtain
$$
 2d(y)+d(z)+d(x_h)$$ $$\leq 3n+1+2r-2p+q-l-k\leq 4n-2-(l-r)-(k-q)-2p <4n-6,
 $$
which contradicts condition $(M)$ and hence the discussion of Subcase 2.1 is completed. \\

\noindent\textbf{Subcase 2.2.} $t\leq h$.

Then $b\leq h-1$ since $h\geq t\geq b+1$.

Assume first that $t=h$. Then  $x_s\rightarrow x_h\rightarrow x_l$. 
By Lemma 4.2, 
$$A(\{x_1,\ldots , x_{h-1}\}\rightarrow \{x_{h+1},\ldots , x_{n-2}\})\not 
=\emptyset.
$$
 Let 
 $x_{i}\rightarrow x_{j}$, where $i\in [1,h-1]$ and $j\in [h+1,n-2]$. 
 From the second equality of (35) it follows that $j\leq k$, i.e., $j\in [h+1,k]$. By
  (36) we have that $i\geq b$, i.e., $i\in [b,h-1]$. Therefore, $C(y,z)=x_1\ldots x_ax_b\ldots x_ix_j\ldots x_kyx_q\ldots x_sx_hx_l\ldots x_{n-2}zx_1$, a contradiction.
 
 Assume next that $t\leq h-1$. From the maximality of $b$ and $t$ it follows that $d^-(x_h,\{x_1,\ldots , x_{b-1}\})=0$. This last equality together with (37) implies that  
 $$
 d(x_h)=d^+(x_h,\{x_1,\ldots , x_{b-1}\})+d(x_h,\{x_b,\ldots , x_k\})+d^-(x_h,\{x_{k+1},\ldots , x_{l-1}\})$$ $$+d(x_h,\{x_{l},\ldots , x_{n-2}\})+d(x_h,\{z\})\leq b-1+2k-2b+l-1-k+2n-2l-2+2$$ $$
 = 2n-l-2+k-b. 
 $$
This  together with (11), $d(z)\leq n-1$ and $r\leq n-3$ implies that
 $$
 2d(y)+d(x_h)+d(z)\leq 2q-2p+2r-2k+4+2n-l-2+k-b+n-1$$ $$\leq 4n-2-(l-r)-(k-q)-(b-q)-2p,
 $$
which contradicts condition $(M)$, since $k-q\geq 0$, $b-q\geq 1$. The discussion of Case 2 is completed. Theorem 1.12 
is proved.    
\end {proof}   


\acknowledgements
\label{sec:ack}
The author would like to thank the anonymous referees  for thoroughly review and many helpful comments and suggestions which improved substantially the rewriting of this paper. We also thank  PhD P. Hakobyan for formatting the manuscript of this paper.

\end{document}